\documentclass[lettersize,journal]{IEEEtran}
\usepackage{amsmath,amsfonts, amssymb}
\usepackage{algorithmic}
\usepackage{algorithm}
\usepackage{array}
\usepackage[caption=false,font=normalsize,labelfont=sf,textfont=sf]{subfig}
\usepackage{textcomp}
\usepackage{stfloats}
\usepackage{url}
\usepackage{verbatim}
\usepackage{graphicx}
\usepackage{cite}
\usepackage{color}
\usepackage{tikz,pgfplots}
\usetikzlibrary{automata, arrows.meta, positioning}

\begin{document}

\title{A Graph-Based Iterative Strategy for Solving the All-Line Transmission Switching Problem}

\author{M. Aguilar-Moreno, S. Pineda,~\IEEEmembership{Senior Member,~IEEE,} J. M. Morales,~\IEEEmembership{Senior Member,~IEEE} 
\thanks{The work of M. Aguilar-Moreno, S. Pineda and J. M. Morales was supported in part by the Spanish Ministry of Science and Innovation (AEI/10.13039/501100011033) through projects PID2020-115460GB-I00 and PID2023-148291NB-I00. M. Aguilar-Moreno, S. Pineda, and J. M. Morales are with the research group OASYS, University of Malaga, Malaga 29071, Spain: marina.aguilar@uma.es; spineda@uma.es; juan.morales@uma.es. Finally, the authors thankfully acknowledge the computer resources, technical expertise, and assistance provided by the SCBI (Supercomputing and Bioinformatics) center of the University of M\'alaga.}
}



\maketitle

\begin{abstract}
The transmission switching problem aims to determine the optimal network topology that minimizes the operating costs of a power system. This problem is typically formulated as a mixed-integer optimization model, which involves big-M constants that lead to weak relaxations and significant computational challenges, particularly when all lines are switchable. In this paper, we propose a two-fold approach: first, using graph theory to derive tighter big-M values by solving a relaxed longest path problem; second, introducing an iterative algorithm that incorporates a heuristic version of the switching problem to efficiently generate low-cost feasible solutions, thereby accelerating the search for optimal solutions in the integer optimization solver. Numerical results on the 118-bus network show that the proposed methodology significantly reduces the computational burden compared to conventional approaches.
\end{abstract}


\section{Introduction}\label{sec:intro}

\IEEEPARstart{T}{ransmission} lines in electricity networks have traditionally been regarded as static infrastructure, with limited control capabilities beyond maintenance periods or outage scenarios. However, the idea of dynamically reconfiguring the grid's topology gained prominence in~\cite{o2005dispatchable} and was later formalized in~\cite{fisher2008optimal}  as the Optimal Transmission Switching (OTS) problem. OTS aims to determine the most efficient configuration of transmission lines in a power system to achieve specific objectives, such as minimizing operating costs, reducing transmission losses, mitigating voltage deviations, or alleviating congestion. Even with the widely adopted direct current (DC) linear approximation of power flow equations, the DC Optimal Transmission Switching (DC-OTS) problem is formulated as a mixed-integer nonlinear program (MINLP), which imposes a significant computational burden on commercial solvers \cite{pineda2024tight}. A common strategy is to reformulate the DC-OTS problem as a mixed-integer linear program (MILP) by using sufficiently large constants. However, the resulting problem has been proven to be NP-hard \cite{kocuk2016cycle}, and finding tight values for these constants remains a challenging task.

One key aspect in determining tight values for these large constants and reducing the computational burden depends on the existence of transmission lines that are permanently connected (non-switchable) and whether the OTS solution must ensure network connectivity by guaranteeing a continuous path of connected lines between all pairs of buses in the system. In the simplest case, a predefined subset of non-switchable lines remains permanently connected, ensuring network connectivity. This scenario is common in practice, as system operators often avoid opening specific lines to ensure system stability. From a mathematical perspective, this case is easier to address for two reasons. First, it involves fewer binary variables, which reduces the search space of the integer model and lowers its computational complexity. Second, as demonstrated in \cite{binato2001new,fattahi2019bound}, the existence of a spanning subtree of connected lines simplifies the computation of bounds on voltage angle differences, which are essential for solving the OTS optimization model as discussed in detail in Section \ref{sec:formulation}.

A more complex scenario arises when all lines in the network are switchable. In this case, the model includes a significantly larger number of binary variables, which expands the search space and increases computational complexity. Additionally, obtaining reliable bounds on voltage angle differences becomes more challenging, further complicating the solution process. Although references \cite{hedman2010co, moulin2010transmission, numan2023role} propose bounds on voltage angles, these bounds are relatively loose, leading to a high computational burden for the OTS problem. Alternatively, other authors propose heuristic algorithms to reduce the number of switchable lines, albeit at the cost of obtaining suboptimal network configurations \cite{liu2012heuristic, barrows2012computationally, fuller2012fast, crozier2022feasible}. Within the category of the all-line OTS problem, the formulation can either include connectivity constraints \cite{li2021connectivity}, ensuring that the resulting network configuration maintains a continuous path between all pairs of buses, or omit these constraints, allowing solutions that result in electrically isolated network segments, commonly known as ``electric islands'' \cite{ostrowski2014transmission}. 

In this paper, we focus on the OTS problem where all transmission lines are switchable, but connectivity constraints are imposed to ensure the network remains fully connected. Specifically, we propose two complementary strategies to reduce the computational burden of this problem. First, we leverage graph theory to efficiently compute tighter bounds on the voltage angle differences. Second, we develop a heuristic procedure to quickly identify high-quality feasible solutions, which are subsequently provided to the mixed-integer algorithm to reduce the search space and the time required to achieve the global optimal solution. The contributions of this paper are thus threefold:
\begin{itemize}
\item[-] We use graph theory to propose an algorithm that computes tighter bounds on voltage angle differences in the OTS problem than those commonly used in the literature. 
\item[-] We develop a heuristic algorithm that efficiently identifies high-quality feasible solutions of the OTS problem, which are iteratively supplied to the mixed-integer algorithm to accelerate its convergence. 
\item[-] We combine these two strategies synergistically and demonstrate their effectiveness through numerical results on a 118-bus network commonly used in OTS studies. 
\end{itemize}

The remainder of this paper is structured as follows. Section \ref{sec:formulation} introduces the OTS formulations and discusses existing approaches for computing big-M values. Section \ref{sec:methodology} provides a detailed explanation of the two proposed strategies for reducing the computational burden of the OTS problem. Section \ref{sec:comparison} outlines the framework used to compare the performance of the proposed methodology with existing approaches. Numerical simulation results using a 118-bus network are presented in Section \ref{sec:numerical_simulations}. Finally, the conclusions are drawn in Section \ref{sec:conclusions}.












\section{Formulation} \label{sec:formulation}

Consider a power network composed of a set of nodes $\mathcal{N}$ and transmission lines $\mathcal{L}$. For simplicity, we assume that each node $n \in \mathcal{N}$ has one generator and one power load. Let $p_n$ and $d_n$ represent the power output of the generator and the power consumption of the load, respectively. Each generator is characterized by its minimum and maximum power output, denoted as $\underline{p}_n$ and $\overline{p}_n$, as well as a marginal production cost $c_n$. The power flow through the transmission line $l=(n,m) \in \mathcal{L}$, connecting nodes $n$ and $m$, is denoted by $f_l$. By convention, $f_l > 0$ indicates a power flow from node $n$ to node $m$, while $f_l < 0$ indicates a flow in the reverse direction. The maximum line capacity is given by $F_l$. Additionally, for any line $l=(n,m) \in \mathcal{L}$, the binary variable $x_l$ indicates its operational status, where $x_l = 1$ means the line is fully operational, and $x_l = 0$ means it is disconnected. Using the DC approximation of the network equations, the power flow $f_l$ through an active line is given by the product of the line's susceptance, $b_l$, and the voltage angle difference between nodes $n$ and $m$, i.e., $\theta_n - \theta_m$. 
%
%
%
With this notation established, the DC-OTS problem can be formulated as follows:
\begin{subequations}\label{eq:OTS_NP}
\begin{IEEEeqnarray}{l}
\min \quad \sum_{n} c_{n} \, p_{n} \label{eq:OTS_NP_obj}\\
\text{subject to}  \\ 
f_l = x_lb_l(\theta_n-\theta_m), \quad \forall l=(n,m) \in \mathcal{L} \label{eq:OTS_NP_Flow_S}\\
p_n - d_n = \sum_{l\in\mathcal{L}(n,\cdot)} f_l - \sum_{l\in\mathcal{L}(\cdot,n)} f_l, \quad \forall n \in \mathcal{N} \label{eq:OTS_NP_PB}\\
\underline{p}_n \leq p_n \leq \overline{p}_n, \quad \forall n \in \mathcal{N} \label{eq:OTS_NP_Plimits}\\
-x_lF_l \leq f_l \leq x_lF_l, \quad \forall l \in \mathcal{L} \label{eq:OTS_NP_Flow_limit_S}\\
\theta_1 = 0 \label{eq:OTS_NP_slack}\\
\sum_{l\in\mathcal{L}(n,\cdot)} y_l + \sum_{l\in\mathcal{L}(\cdot,n)} z_l \geq 1, \quad \forall n \in \mathcal{N} \label{eq:OTS_NP_con_1}\\
u_n - u_m \leq |\mathcal{N}|(1-y_l) - 1, \quad \forall l=(n,m) \in \mathcal{L} \setminus \mathcal{L}_1 \label{eq:OTS_NP_con_2a}\\
u_m - u_n \leq |\mathcal{N}|(1-z_l) - 1, \quad \forall l=(n,m) \in \mathcal{L} \setminus \mathcal{L}_1 \label{eq:OTS_NP_con_2b}\\
y_l + z_l = x_l, \quad \forall l\in\mathcal{L} \label{eq:OTS_NP_con_3}\\
1 \leq u_n \leq |\mathcal{N}|, \quad \forall n \in \mathcal{N} \label{eq:OTS_NP_con_4}\\
u_1 = 1 \label{eq:OTS_NP_con_5}\\
x_l, y_l, z_l \in \{0,1\}, \quad \forall l \in \mathcal{L} \label{eq:OTS_NP_binary}
\end{IEEEeqnarray}
\end{subequations}

The objective function in equation~\eqref{eq:OTS_NP_obj} aims to minimize the total cost of electricity generation. The power flow through the transmission lines is defined in equation~\eqref{eq:OTS_NP_Flow_S}. The nodal power balance is enforced by~\eqref{eq:OTS_NP_PB}, while the power output of each generator is constrained by~\eqref{eq:OTS_NP_Plimits}, ensuring it remains within the limits $[\underline{p}_n, \overline{p}_n]$. Constraint~\eqref{eq:OTS_NP_Flow_limit_S}  imposes limits on the maximum power flow through switchable and non-switchable lines, respectively. The voltage angle of one node is fixed to zero in equation~\eqref{eq:OTS_NP_slack}. Equations~\eqref{eq:OTS_NP_con_1}--\eqref{eq:OTS_NP_con_5} enforce the connectivity of the optimal transmission switching (OTS) solution using the Miller–Tucker–Zemlin (MTZ) formulation~\cite{li2021connectivity}. This formulation relies on a directed graph representation, which requires introducing two binary variables for each line $l$, denoted as $y_l$ and $z_l$. Additionally, the subset $\mathcal{L}_1$ comprises all the lines connected to bus~1. Finally, the binary variables are defined in constraint~\eqref{eq:OTS_NP_binary}.

Problem~\eqref{eq:OTS_NP} is classified as a mixed-integer nonlinear programming (MINLP) problem because of the product term $x_l(\theta_n - \theta_m)$ in~\eqref{eq:OTS_NP_Flow_S}. Constraint~\eqref{eq:OTS_NP_Flow_S} can be linearized by introducing sufficiently large constants $M_{l}$ for each line \cite{hedman2012flexible}, replacing it with  the following two inequalities:
\begin{subequations}\label{eq:linear}
\begin{align}
& f_l \geq b_l(\theta_n-\theta_m)-M_{l}(1-x_{l}), \quad \forall l=(n,m)\in\mathcal{L} \\
& f_l \leq b_{l}(\theta_n-\theta_m)+M_{l}(1-x_{l}), \quad \forall l=(n,m)\in\mathcal{L} 
\end{align}
\end{subequations}
\noindent where the big-M constants $M_l$  are chosen to be upper bounds for $b_{l}|\theta_n - \theta_m|$ when the line $l$ is disconnected ($x_{l} = 0$). With this assumption, the DC-OTS problem can be reformulated as the following mixed-integer linear programming problem:
\begin{subequations} 
\begin{align}
& \min \quad  \sum_{n} c_{n} \, p_{n} \label{eq:ots_mip_of}\\
& \text{subject to} \quad \eqref{eq:OTS_NP_PB}-\eqref{eq:OTS_NP_binary}, \eqref{eq:linear}
\end{align} \label{eq:ots_mip}
\end{subequations}
While model \eqref{eq:ots_mip} can be solved using standard mixed-integer optimization solvers like Gurobi \cite{gurobi}, the selection of appropriate values for the big-M constants $M_{l}$ is critical. As discussed in \cite{fattahi2019bound}, excessively large values of $M_{l}$ can weaken the linear relaxation used in branch-and-bound or branch-and-cut algorithms, leading to increased computational effort. Conversely, smaller values of $M_{l}$ generally result in stronger linear relaxations, which can accelerate convergence. However, as our simulation results indicate in Section \ref{sec:numerical_simulations}, excessively tight big-M values can sometimes lead to longer solution times. This counterintuitive behavior may arise due to numerical stability issues, less effective presolve reductions, or changes in branching decisions that hinder solver performance. Therefore, careful tuning of $M_{l}$ is necessary to balance numerical stability and computational efficiency.


First, we consider the case where there exists a subset of connected non-switchable lines that maintains the network's connectivity. In such case, the authors of \cite{fattahi2019bound} prove that valid big-Ms values can be computed as 
%
%
\begin{equation} \label{eq:bigM_shortest}
    M_{l} = b_{l} \sum_{l'\in {\rm SP}_{l}}\frac{F_{l'}}{b_{l'}}
\end{equation}
\noindent where ${\rm SP}_{l}$ denotes the shortest path between the origin and destination buses of line $l$, which can be efficiently determined using Dijkstra's algorithm~\cite{cormen2022introduction}. Furthermore, these bounds can be further refined using the iterative algorithm proposed in~\cite{pineda2024tight}.

However, as demonstrated in \cite{binato2001new, moulin2010transmission}, when there is no pre-established spanning tree of connected non-switchable lines and all network lines are switchable, these big-M values must be computed as follows:
\begin{equation}
M_{l} = b_{l}  \sum_{l'\in {\rm LP}_{l}}\frac{F_{l'}}{b_{l'}}    
\end{equation}
\noindent where ${\rm LP}_{l}$ denotes the longest path between the origin and destination buses of line $l$. Unlike the shortest path problem, the longest path problem is well-known to be NP-hard \cite{karger1997approximating}, which imposes a prohibitive computational burden. To reduce the computational burden of determining valid big-M values for the Optimal Transmission Switching (OTS) problem when all lines are switchable, efficient methods can be used to compute valid upper bounds on the longest path. For example, a straightforward upper bound is obtained by summing the weights of all edges in the graph. A tighter upper bound can be achieved by adding the weights of the $|\mathcal{N}|-1$ largest edges, which also ensures an upper bound on the longest path. Although these bounds are valid for any demand profile and generation portfolio, they are expected to be quite loose \cite{fattahi2019bound}.

\begin{figure}[]
\centering
\begin{tikzpicture} [node distance = 2cm, auto, scale=0.6, transform shape] 
    \node (n5) [state, initial text = {}] {$n_5$};
    \node (n2) [state, right = of n5] {$n_2$};
    \node (n4) [state, below = of n5] {$n_4$};
    \node (n1) [state, below = of n2] {$n_1$};
    \node (n3) [state, below right = of n2, yshift=0.5cm] {$n_3$};
    \path [thick]
        (n4) edge [color = blue] node [label={[rotate = {90},shift={(0,-0.1)}]$l_{1} (10)$}] {} (n5)
        (n5) edge [color = blue] node {$l_{2} (50)$} (n2)
        (n3) edge node [label={[rotate = {-37},shift={(0,0.05)}]$l_{3} (30)$}]{} (n2)
        (n3) edge [color = blue] node [label={[rotate = {37},shift={(0,0.05)}]$l_{4} (32)$}] {} (n1)
        (n4) edge [color = blue] node {$l_{5} (50)$} (n1)
        (n1) edge node [label={[rotate = {90},shift={(0,-0.1)}]$l_{6} (20)$}] {} (n2);
    \end{tikzpicture}
\caption{Counterexample for methodology proposed in \cite{moulin2010transmission}.}
\label{fig:example_claudia}
\end{figure}
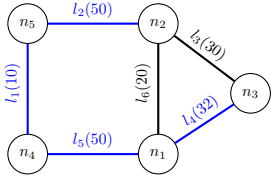

The authors of \cite{moulin2010transmission} propose a greedy algorithm to compute an upper bound of the longest path that, for a given bus ordering, iteratively assigns the edge with the largest weight to each bus that has not yet been connected to a previous bus. However, this approach may yield results smaller than the actual longest path for certain bus orderings, meaning it does not guarantee valid upper bounds. This limitation is illustrated in Figure~\ref{fig:example_claudia}, where the longest path between buses $n_4$ and $n_5$ has a total weight of 162. If the bus ordering is ${n_1, n_2, n_3, n_4, n_5}$, the algorithm assigns edges with weights 50, 50, 32 and 10 to buses $n_1$, $n_2$, $n_3$, and $n_4$, respectively, while no edge is assigned to $n_5$. This results in a total weight of 142, which is smaller than the longest path, demonstrating that the method proposed in \cite{moulin2010transmission} may fail to provide valid upper bounds for the longest path within a graph.

\section{Methodology} \label{sec:methodology}

In this paper, we propose reducing the computational burden of the OTS problem \eqref{eq:ots_mip} through two complementary and synergistic strategies that achieve the best results when combined. First, in Section \ref{sec:tighter_bigm}, we leverage graph theory to simplify the network and employ a relaxed formulation of the longest-path problem, which can be solved efficiently, to derive tight values for $M_l$. Second, in Section \ref{sec:iterative}, we introduce an iterative procedure designed to accelerate the solution of the OTS problem \eqref{eq:ots_mip} by identifying high-quality feasible solutions. These solutions effectively reduce the size of the binary search tree required to solve the mixed-integer problem \eqref{eq:ots_mip}. While the methodology for finding new feasible solutions is heuristic, the iterative procedure guarantees convergence to the global optimal solution. Although the iterative algorithm can be implemented using any big-M values, the lowest computational burden is achieved when the tighter big-Ms derived in Section \ref{sec:tighter_bigm} are utilized.

\subsection{Tighter big-Ms} \label{sec:tighter_bigm}

In this section, we begin by representing the electricity network as a graph, which provides a structured framework for analyzing its topology using graph theory. This representation models the connections between nodes and lines as vertices and edges, respectively, enabling systematic simplifications of the graph to reduce the values of \(M_l\). Following these simplifications, we formulate and solve a relaxed version of the longest-path problem on the simplified graph, which further refines the values of \(M_l\). 

Specifically, we define a weighted multigraph \( G = (\mathcal{N}, \mathcal{L}, \phi, w) \), where \( \mathcal{N} \) is the set of nodes corresponding to the buses in the network, \( \mathcal{L} \) is the set of edges representing the transmission lines, \( \phi: \mathcal{L} \to \{\{n, m\} \mid n, m \in \mathcal{N}, n \neq m\} \) is the incidence function that maps each line to the pair of nodes it connects, and \( w: \mathcal{L} \to \mathbb{R} \) assigns a weight to each line. A node \( n \in \mathcal{N} \) is a neighbor of node \( m \in \mathcal{N} \) if there exists a line \( l \in \mathcal{L} \) such that \( \phi(l) = \{ n, m \} \). The degree of a node \( n \in \mathcal{N} \), denoted \( \deg(n) \), is the number of neighbors of \( n \). A node \( n \in \mathcal{N} \) is called a leaf node if its degree is 1, meaning that it is connected to exactly one other node. A line \( l \in \mathcal{L} \) is called a leaf line if one of its nodes is a leaf node. A node \( n \in \mathcal{N} \) is called a cut node or articulation point if its removal, along with all lines incident to it, disconnects the graph. Two lines \( l_1, l_2 \in \mathcal{L} \) are said to be parallel lines if they connect the same pair of nodes, that is, \( \phi(l_1) = \phi(l_2) \).

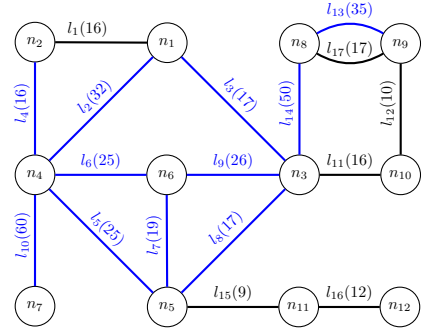
\begin{figure}
    \centering    
    \begin{tikzpicture} [node distance = 3.25cm, auto, scale=0.6, transform shape] 
        \node (n1) [state, initial text = {}] {$n_1$};
        \node (n4) [state, below left = of n1] {$n_4$};
        \node (n2) [state, above  = of n4, yshift=-1.2cm] {$n_2$};
        \node (n3) [state, below right = of n1] {$n_3$};
        \node (n5) [state, below right = of n4] {$n_5$};
        \node (n6) [state, left = of n3, xshift =1.2cm] {$n_6$};
        \node (n7) [state, below = of n4, yshift =1.23cm] {$n_7$};
        \node (n11) [state, right = of n5, xshift =-1.23cm] {$n_{11}$};
        \node (n12) [state, right = of n11,xshift =-1.9cm] {$n_{12}$};
        \node (n8) [state, above = of n3, yshift=-1.23cm] {$n_8$};
        \node (n9) [state, right = of n8, xshift =-1.9cm] {$n_9$};
        \node (n10) [state, right = of n3, xshift =-1.9cm] {$n_{10}$};
        \path [thick]
            (n1) edge node [label={[shift={(-0.3,0)}]$l_{1} (16)$}] {} (n2)
            (n3) edge [color = blue] node [label={[rotate = {-45},shift={(0,0.1)}]$l_{3} (17)$}] {} (n1)
            (n4) edge [color = blue] node [label={[rotate = {45},shift={(0,-0.3)}]$l_{2} (32)$}] {} (n1)
            (n4) edge [color = blue] node [label={[rotate = {90},shift={(0,-0.15)}]$l_{4} (16)$}] {} (n2)
            (n4) edge [color = blue] node [label={[rotate = {-45},shift={(0,-0.3)}]$l_{5} (25)$}] {} (n5)
            (n3) edge [color = blue] node [label={[rotate = {45},shift={(0,0.1)}]$l_{8} (17)$}] {} (n5)
            (n4) edge [color = blue] node {$l_{6} (25)$} (n6)
            (n3) edge [color = blue] node [label={[shift={(0,0)}]$l_{9} (26)$}]{} (n6)
            (n5) edge [color = blue] node [label={[rotate = {90},shift={(0,-0.15)}]$l_{7} (19)$}] {} (n6)
            (n4) edge[color = blue]  node [label={[rotate = {90},shift={(0,0.1)}]$l_{10} (60)$}] {}(n7)
            (n5) edge node {$l_{15} (9)$} (n11)
            (n11) edge node {$l_{16} (12)$} (n12)
            (n3) edge [color = blue] node [label={[rotate = {90},shift={(0,-0.15)}]$l_{14} (50)$}] {} (n8)
            (n3) edge node {$l_{11} (16)$} (n10)
            (n10) edge node [label={[rotate = {90},shift={(0,-0.15)}]$l_{12} (10)$}] {} (n9)
            (n8) edge [color = blue, bend left] node {$l_{13} (35)$} (n9)
            (n8) edge [bend right] node {$l_{17} (17)$} (n9);
        \end{tikzpicture}
    \caption{Illustrative example for determining \(M_{l_1}\) using the naive approach.}
    \label{fig:naive}
\end{figure}

To illustrate the simplifications we propose, we refer to Figure \ref{fig:naive}, which depicts a graph containing a set of nodes \( n_1, n_2, \dots , n_{12} \), and a set of lines \( l_1, l_2, \dots, l_{17} \) with their corresponding weights in parentheses. Vertices $n_7$ and $n_{12}$ are leaf nodes, and edges $l_{10}$ and $l_{16}$ are leaf lines. Additionally, node $n_3$ is a cut vertex, and $l_{13}$ and $l_{17}$ are parallel lines.

To illustrate the graph simplifications we propose, we particularize it to compute the values $M_l$ associated with the line $l_1$ that connects buses $n_1$ and $n_2$. Since the graph contains 12 vertices, the naive approach, which adds the largest \( |\mathcal{N}| - 1 \) weights marked in blue, results in the following big-M value:
\[ M_{l_1} = b_{l_1} \cdot (60 + 50 + 35 + \dots + 17 + 17 + 16) = b_{l_1}\cdot 322 \]

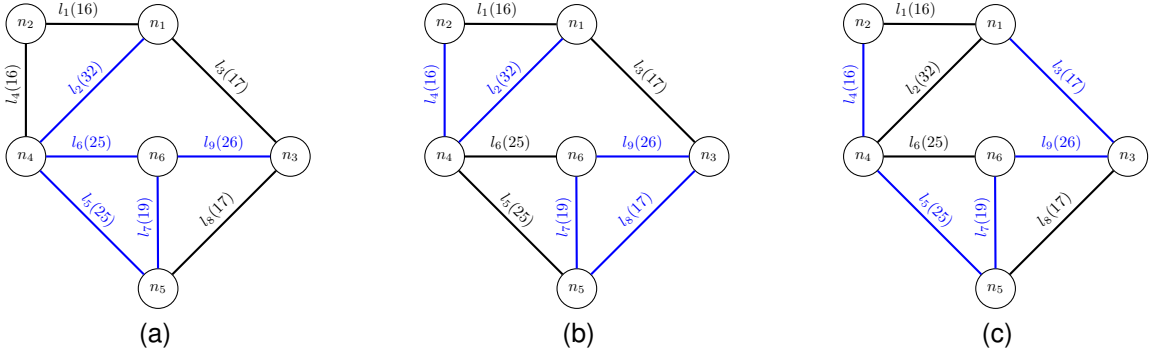
\begin{figure*}[!t]
\centering
\subfloat[]{    \begin{tikzpicture} [node distance = 3.25cm, auto, scale=0.6, transform shape] 
    \node (n1) [state, initial text = {}] {$n_1$};
    \node (n4) [state, below left = of n1] {$n_4$};
    \node (n2) [state, above  = of n4, yshift=-1.2cm] {$n_2$};
    \node (n3) [state, below right = of n1] {$n_3$};
    \node (n5) [state, below right = of n4] {$n_5$};
    \node (n6) [state, left = of n3, xshift =1.2cm] {$n_6$};
    \path [thick]
        (n1) edge node [label={[shift={(-0.3,0)}]$l_{1} (16)$}] {} (n2)
        (n3) edge node [label={[rotate = {-45},shift={(0,0.1)}]$l_{3} (17)$}] {} (n1)
        (n4) edge [color = blue] node [label={[rotate = {45},shift={(0,-0.3)}]$l_{2} (32)$}] {} (n1)
        (n4) edge node [label={[rotate = {90},shift={(0,-0.15)}]$l_{4} (16)$}] {} (n2) 
        (n4) edge [color = blue] node [label={[rotate = {-45},shift={(0,-0.3)}]$l_{5} (25)$}] {} (n5)
        (n3) edge node [label={[rotate = {45},shift={(0,0.1)}]$l_{8} (17)$}] {} (n5)
        (n4) edge [color = blue] node {$l_{6} (25)$} (n6)
        (n3) edge [color = blue] node [label={[shift={(0,0)}]$l_{9} (26)$}]{} (n6)
        (n5) edge [color = blue] node [label={[rotate = {90},shift={(0,-0.15)}]$l_{7} (19)$}] {} (n6); 
    \end{tikzpicture}
\label{fig:naive_simpl}}
\hfil
\subfloat[]{\begin{tikzpicture} [node distance = 3.25cm, auto, scale=0.6, transform shape] 
    \node (n1) [state, initial text = {}] {$n_1$};
    \node (n4) [state, below left = of n1] {$n_4$};
    \node (n2) [state, above  = of n4, yshift=-1.2cm] {$n_2$};
    \node (n3) [state, below right = of n1] {$n_3$};
    \node (n5) [state, below right = of n4] {$n_5$};
    \node (n6) [state, left = of n3, xshift =1.2cm] {$n_6$};
    \path [thick]
        (n1) edge node [label={[shift={(-0.3,0)}]$l_{1} (16)$}] {} (n2)
        (n3) edge node [label={[rotate = {-45},shift={(0,0.1)}]$l_{3} (17)$}] {} (n1)
        (n4) edge [color = blue] node [label={[rotate = {45},shift={(0,-0.3)}]$l_{2} (32)$}] {} (n1)
        (n4) edge [color = blue] node [label={[rotate = {90},shift={(0,-0.15)}]$l_{4} (16)$}] {} (n2) 
        (n4) edge node [label={[rotate = {-45},shift={(0,-0.3)}]$l_{5} (25)$}] {} (n5)
        (n3) edge [color = blue] node [label={[rotate = {45},shift={(0,0.1)}]$l_{8} (17)$}] {} (n5)
        (n4) edge node {$l_{6} (25)$} (n6)
        (n3) edge [color = blue] node [label={[shift={(0,0)}]$l_{9} (26)$}]{} (n6)
        (n5) edge [color = blue] node [label={[rotate = {90},shift={(0,-0.15)}]$l_{7} (19)$}] {} (n6);
    \end{tikzpicture}
\label{fig:LP_relaxed}}
\hfil
\subfloat[]{\begin{tikzpicture} [node distance = 3.25cm, auto, scale=0.6, transform shape] 
    \node (n1) [state, initial text = {}] {$n_1$};
    \node (n4) [state, below left = of n1] {$n_4$};
    \node (n2) [state, above  = of n4, yshift=-1.2cm] {$n_2$};
    \node (n3) [state, below right = of n1] {$n_3$};
    \node (n5) [state, below right = of n4] {$n_5$};
    \node (n6) [state, left = of n3, xshift =1.2cm] {$n_6$};
    \path [thick]
        (n1) edge node [label={[shift={(-0.3,0)}]$l_{1} (16)$}] {} (n2)
        (n3) edge [color = blue] node [label={[rotate = {-45},shift={(0,0.1)}]$l_{3} (17)$}] {} (n1)
        (n4) edge node [label={[rotate = {45},shift={(0,-0.3)}]$l_{2} (32)$}] {} (n1)
        (n4) edge [color = blue] node [label={[rotate = {90},shift={(0,-0.15)}]$l_{4} (16)$}] {} (n2) 
        (n4) edge [color = blue] node [label={[rotate = {-45},shift={(0,-0.3)}]$l_{5} (25)$}] {} (n5)
        (n3) edge node [label={[rotate = {45},shift={(0,0.1)}]$l_{8} (17)$}] {} (n5)
        (n4) edge node {$l_{6} (25)$} (n6)
        (n3) edge [color = blue] node [label={[shift={(0,0)}]$l_{9} (26)$}]{} (n6)
        (n5) edge [color = blue] node [label={[rotate = {90},shift={(0,-0.15)}]$l_{7} (19)$}] {} (n6);  
    \end{tikzpicture}
\label{fig:LP}}
\caption{Illustrative example for determining \(M_{l_1}\) using different approaches.}
\label{fig:comparison_ml1}
\end{figure*}


Upon examining the solution obtained from this naive approach, it becomes evident that, although the leaf edge \( l_{10} \) has a large weight, it cannot be part of the longest path between the source node \( n_1 \) and the sink node \( n_2 \). This is because the longest path between \( n_1 \) and \( n_2 \) must be simple (i.e., it cannot repeat vertices), and as such, it cannot include leaf edges, which would terminate the path prematurely. Therefore, a natural graph simplification is to iteratively remove all leaf edges until none remain. Specifically, in the graph shown in Figure \ref{fig:naive_simpl}, edges $l_{10}$, $l_{15}$, and $l_{16}$ are removed, and the big-M values associated with these edges are set to 0, indicating that, in the OTS solution, these lines are always connected.

A second observation concerns the existence of cut node $n_3$, which divides the graph into two subgraphs. Although edges $l_{13}$ and $l_{14}$ have high weights, they cannot be part of the longest path from $n_1$ to $n_2$, as they belong to the right-hand subgraph, and the cut vertex $n_3$ can only be visited once in the longest path solution. This implies that when computing the big-M values for an edge in one subgraph, none of the edges from the other subgraphs should be included in the longest path. In other words, if a graph contains \( |\mathcal{N}_C| \) cut edges, the graph can be divided into \( |\mathcal{N}_C| + 1 \) subgraphs, and the longest paths must be computed independently for each subgraph.

If the leaf edges \( l_{10} \), \( l_{15} \), and \( l_{16} \) are removed, and we consider only the left-hand subgraph, which contains 9 edges and 6 nodes, the new value for the big-M parameter is determined by the five edges with the highest weights as:
\[ M_{l_1} = b_{l_1} \cdot (32+26+ 25+ 25+ 19) = b_{l_1}\cdot 127 \]
\noindent which is 60\% lower than the original value. The selected lines after these simplifications are shown in Figure \ref{fig:naive_simpl}. However, the selected lines do not fully satisfy the longest path requirements, as each node in the path should be connected to exactly two lines, one for entering and one for leaving. In the current selection of lines of Figure \ref{fig:naive_simpl}, some nodes (\(n_4\) and \(n_6\)) are connected to more than two lines, which violates this requirement. Therefore, this bound can be improved by excluding solutions where any node is connected to more than two selected lines. To achieve this, we propose determining the big-M values by solving a relaxed version of the longest-path problem on the simplified graph, obtained by removing leaf edges and partitioning the graph. Below is the formulation of the longest-path problem without sub-tour elimination constraints, which simplifies the computational complexity by excluding these constraints, as discussed in \cite{taccari2016integer}.
\begin{subequations}\label{eq:LP}
\begin{IEEEeqnarray}{l}
\max_{y_{l},z_{l}} \quad \sum_{l} \frac{F_{l}}{b_{l}} \, (y_{l} + z_{l}) \label{eq:LP_obj}\\
\text{subject to} \\  
\sum_{l\in\mathcal{L}(\cdot,n)} (y_{l} + z_{l}) + \beta_n = \sum_{l\in\mathcal{L}(n,\cdot)} (y_{l} + z_{l}), \quad \forall n \in \mathcal{N} \label{eq:LP_Flow}\\
\sum_{l\in\mathcal{L}(n,\cdot)} (y_{l} + z_{l}) \leq 1, \quad \forall n \in \mathcal{N} \mid n \neq sink \label{eq:LP_out} \\
\sum_{l\in\mathcal{L}(n,\cdot)} (y_{l} + z_{l}) = 0, \quad \forall n \in \mathcal{N} \mid n = sink \label{eq:LP_out_sink} \\
y_{l} + z_{l} \leq 1, \quad \forall l \in \mathcal{L} \label{eq:LP_binary_sum} \\
y_{l},z_{l} \in \{0,1\}, \quad \forall l \in \mathcal{L} \label{eq:LP_binary} 
\end{IEEEeqnarray}
\end{subequations}
%
\noindent where parameter $\beta_n$ is equal to 1 for the sources node, to -1 for the sink node, and 0 otherwise. If problem \eqref{eq:LP} is solved for line $l_1$ with the source node $n_1$ and sink node $n_2$, we obtain the solution shown in Figure \ref{fig:LP_relaxed}. Note that in this solution, none of the vertices is connected to more than two selected edges. However, the solution is still not a valid path, as it contains subtours. Therefore, the obtained solution provides an upper bound for the longest path, but it represents a tighter bound than the one shown in Figure \ref{fig:naive_simpl}. In fact, the value of the $M_{l_1}$ for this approach is $b_{l_1} \cdot 110$, which is 66\% lower than the original value. For completeness, Figure \ref{fig:LP} illustrates the solution to the computationally more intensive longest path problem, yielding \(M_{l_1} = b_{l_1} \cdot 103\), which is slightly tighter than the relaxed solution shown in Figure \ref{fig:LP_relaxed}.

Therefore, through this illustrative example, we have demonstrated how the big-M values can be tightened using graph properties. In summary, the procedure we propose proceeds as follows: First, we iteratively remove all leaf edges and assign big-M values of 0 to those edges. Second, we split the graph into subgraphs if cut vertices exist. Third, we solve the optimization model \eqref{eq:LP} for each line $l$, considering only the weights of the subgraph to which this line belongs. Using these lower big-M values directly to solve the OTS problem \eqref{eq:ots_mip} results in a tighter relaxations that are expected to reduce the computational burden.  

\subsection{Iterative algorithm} \label{sec:iterative}

In this section, we propose a simple and fast heuristic methodology to efficiently find feasible solutions that may potentially improve the incumbent solution found by the mixed-integer solver. For given big-M values, the proposed approach begins by solving the optimization model \eqref{eq:ots_mip} using a mixed-integer solver. After a pre-established time limit (e.g., 100 seconds), the solver is stopped, and all the feasible solutions obtained so far are retrieved. Among them, we select a subset of feasible solutions $\mathcal{K}$, indexed by \( k \), with the lowest objective function values. The solutions in this subset are characterized by the connection status of all lines denoted by $x_{kl}, \forall k \in \mathcal{K}, \forall l \in \mathcal{L}$. All these solutions are feasible in model \eqref{eq:ots_mip} and ensure network connectivity by including a path of connected lines between any two buses.

The heuristic procedure we propose to find new feasible solutions involves solving a simplified version of the optimization model \eqref{eq:ots_mip}, where binary variables are fixed to 1 for lines consistently connected across all solutions in the subset \( \mathcal{K} \). Additionally, for each solution in \( \mathcal{K} \), big-M values are computed using a shortest-path approach, and the largest big-M value for each line across all solutions in the subset is selected. By applying these two simplifications, we obtain the following modified formulation of the OTS problem:
\begin{subequations}\label{eq:OTS_NP_heur}
\begin{IEEEeqnarray}{l}
\min \quad \sum_{n} c_{n} \, p_{n} \label{eq:OTS_NP_heur_obj}\\
\text{subject to} \\ 
\eqref{eq:OTS_NP_PB}-\eqref{eq:OTS_NP_binary} \\
f_l \geq b_l(\theta_n-\theta_m)-(1-x_{l})\max_k(M_{lk}), \forall l\in\mathcal{L} \label{eq:OTS_NP_heur_bigm1}\\
f_l \leq b_{l}(\theta_n-\theta_m)+(1-x_{l})\max_k(M_{lk}), \forall l\in\mathcal{L} \label{eq:OTS_NP_heur_bigm2}\\
\left\lfloor \frac{\sum_k x_{kl}}{|\mathcal{K}|} \right\rfloor \leq x_l \leq 1, \quad \forall l \in \mathcal{L} \label{eq:OTS_NP_heur_fix_binary}
\end{IEEEeqnarray}
\end{subequations}


Since formulation \eqref{eq:OTS_NP_heur} includes fewer binary variables and tighter big-M values compared to formulation \eqref{eq:ots_mip}, it is computationally easier to solve. However, the optimal solution of \eqref{eq:OTS_NP_heur} is not guaranteed to be the optimal solution of the OTS formulation \eqref{eq:ots_mip}, as the constraints \eqref{eq:OTS_NP_heur_bigm1}–\eqref{eq:OTS_NP_heur_fix_binary} are determined heuristically. Nevertheless, it is important to note that all solutions in \( \mathcal{K} \) satisfy equations \eqref{eq:OTS_NP_heur_bigm1}–\eqref{eq:OTS_NP_heur_fix_binary}, as they are feasible solutions of model \eqref{eq:ots_mip}, and they also satisfy constraints \eqref{eq:OTS_NP_PB}-\eqref{eq:OTS_NP_binary} by definition, ensuring that all solutions in \( \mathcal{K} \) are feasible for model \eqref{eq:OTS_NP_heur}. Furthermore, the optimal solution of model \eqref{eq:OTS_NP_heur} is always feasible for the original OTS model \eqref{eq:ots_mip}, as it shares the same objective function but has a larger feasible region. 

Based on these properties, we propose an iterative procedure that alternates between solving the original OTS model \eqref{eq:ots_mip} and its simplified, heuristic version \eqref{eq:OTS_NP_heur}. For a pre-defined number of iterations, indexed by \( i \), we define \( T_i^O \) as the maximum computational time allocated to solve the original OTS problem \eqref{eq:ots_mip}, and \( T_i^H \) as the maximum computational time allocated to solve the heuristic version \eqref{eq:OTS_NP_heur} during the \( i \)-th iteration. Additionally, \( |\mathcal{K}_i| \) represents the size of the subset of suboptimal solutions selected in the \( i \)-th iteration. At each iteration \( i \), the optimization problem \eqref{eq:ots_mip} is first solved within a maximum computational time of \( T_i^O \). If the optimal solution is found within this time, the algorithm terminates. Otherwise, the \( |\mathcal{K}_i| \) feasible solutions with the lowest objective function values are selected. The big-M values for each solution in subset \(\mathcal{K}_i\) are then determined using equation \eqref{eq:bigM_shortest}, based on the shortest path formulation. Next, the simplified model \eqref{eq:OTS_NP_heur} is solved using the best feasible solution from \( \mathcal{K}_i \) as the starting point. The solutin of the heuristic model \eqref{eq:OTS_NP_heur} within the maximum computational time of \( T_i^H \) is then used to initialize the next iteration of solving \eqref{eq:ots_mip}. A pseudo-code of this iterative methodology is provided in Algorithm \ref{alg:proposed}.

\begin{algorithm}[H]
  \caption{Pseudo-code of proposed approach} \label{alg:proposed}
\textbf{Input:} Big-M values $\mathbf{M}^O$ computed by BN/SR/LP, maximum time $T$, number of iterations $|\mathcal{I}|$, maximum times per iteration $T^O_i, T^H_i \forall i \in \mathcal{I}$, and the size of solution subsets $|\mathcal{K}_i| \forall i \in \mathcal{I}$.

\textbf{Output:} Solution of OTS model \eqref{eq:ots_mip}.

 \begin{algorithmic}[1]
    \STATE Set $i=1$ and $\mathbf{x}^O_i=\mathbf{1}$ (all lines connected). 
    \STATE Solve \eqref{eq:ots_mip} with $\mathbf{x}^O_i$ as initial solution, big-M values $\mathbf{M}^O$ and maximum time  $T^O_i$. If solved, terminate. Otherwise, go to step 3.
    \STATE Set $\mathcal{K}_i$ with the best feasible solutions of \eqref{eq:ots_mip}, $\mathbf{x}^H_i$ as the current incumbent.
    \STATE Use subset $\mathcal{K}_i$ to compute  big-M values $\mathbf{M}^H_i$ using \eqref{eq:bigM_shortest}.
    \STATE Solve \eqref{eq:OTS_NP_heur} with $\mathbf{x}^H_i$ as initial solution, big-M values $\mathbf{M}^H_i$, and  maximum time  $T^H_i$. Set $t^H_i$ as the time to solve \eqref{eq:OTS_NP_heur}.      
    \STATE Set $i = i+1$ and $\mathbf{x}^O_i$ to the incumbent of \eqref{eq:OTS_NP_heur}. If $i\leq |\mathcal{I}|$ go to step 2, otherwise go to step 7.    
    \STATE Compute remaining available time as $T - \sum_i T^O_i + t^H_i$.
    \STATE Solve \eqref{eq:ots_mip} with $\mathbf{x}^O_i$ as initial solution, big-M values $\mathbf{M}^O$ and maximum time computed in step 7.
\end{algorithmic}
\end{algorithm}

The proposed iterative algorithm guarantees the global optimal solution for the OTS problem, as it ensures that the optimal solution of the original model \eqref{eq:ots_mip} is eventually obtained. It includes a heuristic step that may occasionally yield a solution coinciding with the best available solution in the subset \( \mathcal{K} \), indicating an unsuccessful iteration. However, the impact on computational time is minimal, as the next iteration starts from the best feasible solution. If the heuristic solution is better than the current best, it has the potential to significantly reduce the binary tree search needed to solve the original OTS model. This iterative approach not only accelerates the solution process but can also be applied using any given big-M values, regardless of whether they are obtained through the graph simplifications discussed in Section \ref{sec:tighter_bigm}.




\section{Comparison} \label{sec:comparison}

In this section, we describe the different methodologies compared in this study to solve the OTS problem when all network lines are switchable. Specifically, we consider three methods for computing the big-M values. The benchmark approach (BN) simply adds the largest \( |N| - 1 \) weights. While this method is extremely fast, it produces very loose big-M values. The longest-path approach (LP) solves an NP-hard problem for each line, resulting in a high computational burden but yielding the tightest big-M values. Finally, we propose the SR approach, which combines graph simplifications (S) with a relaxation (R) of the longest-path problem. This approach is expected to require low computational effort while achieving tighter big-M values, closer to those obtained by the LP approach.  Regardless of the big-M values used, the optimization model \eqref{eq:ots_mip} can be solved either in a single step (SS) or using the iterative procedure (IT) described in Section \ref{sec:iterative}. By combining the three methods for computing big-Ms with the two approaches for solving model \eqref{eq:ots_mip}, this study evaluates six different methodologies for solving the OTS problem, namely, BN-SS, BN-IT, LP-SS, LP-IT, SR-SS, and SR-IT. We compare the six approaches based on their computational burden, specifically measuring the time required to solve the OTS problem and the MIP gap when the problem is not solved to optimality within the given time limit.

\section{Numerical simulations}\label{sec:numerical_simulations}

In this section, we present the numerical simulations corresponding to the different approaches for solving the OTS problem on the commonly used 118-bus network \cite{blumsack2006network}, which includes 186 lines. This network is sufficiently large to challenge current algorithms while remaining computationally manageable. Moreover, the 118-bus system is a paradigmatic benchmark widely employed in numerous studies on optimal transmission switching in the technical literature \cite{fisher2008optimal, kocuk2016cycle, fattahi2019bound, fuller2012fast, crozier2022feasible, hinneck2022optimal, johnson2020knearest, yang2019line, dey2022node}. To ensure a diverse range of problem complexities, we run all the approaches using 100 demand samples generated by randomly sampling the nodal demand from a uniform distribution within the interval $[0.9\widehat{d}_n, 1.1\widehat{d}_n]$, where $\widehat{d}_n$ represents the baseline demand. In the iterative methodology, the maximum time for solving the original OTS problem increases progressively with each iteration (30, 60, 120, 300, and 600 seconds), while the time for the simplified heuristic problem remains fixed at 30 seconds. At each iteration, the 10 best feasible solutions based on the objective function are selected. All optimization problems are solved using GUROBI 10.0.3 \cite{gurobi} on a Linux-based server with CPUs clocked at 2.6 GHz, utilizing a single thread and 8 GB of RAM. The optimality gap is set to 0.01\%, with a 1-hour time limit, integrality focus set to 1, and all other parameters set to their default values.

In Section \ref{sec:intro}, we explain that while solving the OTS problem with all lines switchable is more computationally challenging, it offers significant operational cost reductions. Comparing the full OTS model to a simplified model with a predefined spanning tree shows an average cost saving of 10.35\%, with savings ranging from 3.35\% to 16.01\%, highlighting the benefits of the more flexible, all-line switchable approach.

The numerical results of this section are presented as follows. First, we analyze the big-M values obtained by the BN, LP, and SR approaches described in Section \ref{sec:tighter_bigm}, along with the computational time required to compute them. Second, for a specific demand scenario, we provide detailed results of the iterative approach described in Section \ref{sec:iterative}. Finally, we compare the computational burden of the six approaches described in Section \ref{sec:comparison}, averaged over the 100 instances.

We begin the analysis of the computational results by comparing the big-M values obtained using the methodology proposed in Section \ref{sec:tighter_bigm} with those computed by the BN and LP approaches. Since the big-Ms determined through the longest-path (LP) approach represent the tightest valid upper bounds for the voltage angle differences, we define the following two comparison factors to evaluate the big-Ms of the BN and SR approaches relative to those of the LP approach:
\[ \lambda^{BN}_l = \frac{M^{BN}_l}{M^{LP}_l} \qquad \qquad \lambda^{SR}_l = \frac{M^{SR}_l}{M^{LP}_l} \]
As previously discussed, it follows that $\lambda^{BN}_l \geq \lambda^{SR}_l \geq 1$, indicating that the big-M values derived using the BN approach are at least as large as those from the SR approach, which, in turn, are no smaller than the tightest bounds established by the LP approach. The minimum, average, and maximum values of \(\lambda^{BN}_l\) are 1.88, 3.23, and 24.26, respectively, while for \(\lambda^{SR}_l\), the corresponding values are 1.00, 1.12, and 1.45. It is clear that the proposed approach, which utilizes graph theory to simplify the network and solves a relaxed version of the longest-path problem, yields big-M values that are significantly lower than those obtained using the naive strategy. Furthermore, these values are only slightly larger than those derived from the exact longest-path algorithm. 

Another noteworthy aspect of the proposed SR approach for computing big-M values is its significantly lower computational time compared to the exact longest-path problem. On average, the SR method computes the big-M values in just 0.007 seconds, whereas the exact algorithm requires 4534 seconds and, in some cases, reaches a maximum runtime of 24 hours. Furthermore, for certain lines in the network, the exact longest-path problem terminated prematurely due to memory constraints rather than the time limit, underscoring the considerable computational burden of this approach.

Next, we analyze the impact of big-M values on the computational burden of the OTS problem by comparing the results of the single-step methods BN-SS, LP-SS, and the proposed SR-SS. As previously discussed, selecting the tightest big-M values provided by the longest path approach can be counterproductive, as it may lead to numerical instability, less effective presolve reductions, or changes in branching decisions that degrade solver performance. To further investigate this issue, we also include computational results for the longest path approach with progressively increased big-M values, ranging from 5\% to 40\%, with the specific increase indicated in brackets.

Table \ref{tab:summary_results_ss} presents the average computational time over 100 instances, along with the average and maximum optimality gap and the total number of unsolved instances after one hour for each single-step approach. First, utilizing the big-M values obtained through the graph simplifications proposed in Section \ref{sec:tighter_bigm} results in solving 7 more instances and reducing the average computational time by 12.3\% compared to the benchmark approach. Second, the use of the smallest big-M values derived from the longest path approach does not yield the best computational performance, leading to 55 unsolved instances within one hour and a computational time larger than those of BN-SS and SR-SS. Third, uniformly increasing the big-M values computed via the longest path method may or may not result in computational savings. For example, a 35\% increase in big-M values leads to a lower computational time than the proposed approach, but for any other increase factor, the computational time is higher. These findings suggest that computing big-M values using the longest path approach is impractical, as it requires substantial computational effort while failing to guarantee the lowest computational time for solving the OTS problem in practice.

\begin{table}[]
    \centering
    \begin{tabular}{lcccc}
    \hline
    Method  & Av. time (s) & Max Gap & Av. Gap & \# unsolved \\
    \hline
    BN-SS & 1893.57   & 2.91\%     & 0.18\%      & 47  \\
    SR-SS & 1660.14   & 7.27\%     & 0.18\%      & 40  \\
    LP-SS & 2133.58   & 2.72\%     & 0.19\%      & 55  \\
    LP-SS(+5\%) & 2042.22 & 3.79\%  & 0.27\%      & 54  \\
    LP-SS(+10\%) & 1828.99 & 2.42\%  & 0.09\%      & 47  \\
    LP-SS(+15\%) & 2233.67 & 3.67\%  & 0.25\%      & 57  \\
    LP-SS(+20\%) & 2109.35 & 3.53\%  & 0.31\%      & 55  \\
    LP-SS(+25\%) & 1814.92 & 2.72\%  & 0.13\%      & 45  \\
    LP-SS(+30\%) & 1883.47 & 2.56\%  & 0.15\%      & 48  \\
    LP-SS(+35\%) & 1614.48 & 3.56\%  & 0.17\%      & 39  \\
    LP-SS(+40\%) & 2158.45 & 3.21\%  & 0.25\%      & 56  \\   
    \hline
    \end{tabular}
    \caption{Summary of computational results for single step approaches.}
    \label{tab:summary_results_ss}
\end{table}

In the second part of this section, we analyze the performance of the proposed iterative methodology. To illustrate its effectiveness, we begin by examining a specific demand scenario as an example. Table \ref{tab:iterative_procedure} presents the results corresponding to this demand scenario, detailing the objective values and gaps for both the original OTS problem \eqref{eq:ots_mip} and its heuristic version \eqref{eq:OTS_NP_heur} over five iterations. The final row of the table corresponds to step 8 of Algorithm \ref{alg:proposed}. A missing gap value indicates that the respective model is solved to optimality within the allocated time limit. If the original OTS problem \eqref{eq:ots_mip} is solved to optimality, the iterative procedure terminates, and the optimal solution is obtained. Otherwise, the process continues. In the heuristic step, even if model \eqref{eq:OTS_NP_heur} is not solved to optimality within the time limit, it can still provide a better incumbent solution than the current one. Improvements in the objective function during the heuristic step have been highlighted in bold for clarity.

Results in Table \ref{tab:iterative_procedure} demonstrate the effectiveness of the proposed iterative methodology, which alternates solving the OTS problem with the heuristic model to efficiently identify better feasible solutions. This approach performs well across all three big-M value strategies used in the simulations. Specifically, the heuristic step successfully improves the incumbent solution in three iterations for the BN approach, and in two iterations for the SR and LP approaches. Thanks to these improvements, the iterative methodology achieves significantly faster solution times for the OTS problem under this demand scenario. The total solution times are reduced to 1595, 243, and 163 seconds for the BN, SR, and LP approaches, respectively, compared to 3600, 1336, and 331 seconds required by their single-step counterparts. These results clearly show that the iterative approach is highly effective in reducing computational times while still ensuring the optimal solution of the OTS problem.

\begin{table}[]
    \centering
    \setlength{\tabcolsep}{4pt}
    \begin{tabular}{ccccccccc}
    \hline
         &            & \multicolumn{2}{c}{Benchmark (BN)}                     & \multicolumn{2}{c}{Proposed (SR)}                       & \multicolumn{2}{c}{LongestPath (LP)}                  \\
\hline
Iteration & Model & ObjVal & Gap &     ObjVal& Gap                     &    ObjVal& Gap                  \\
\hline
$i=1$&\eqref{eq:ots_mip}           & 1496.4& 5.88\%       & 1468.6 &4.07\%      & 1505.9 & 6.47\%     \\
$i=1$&\eqref{eq:OTS_NP_heur}       & \textbf{1468.3}& 0.79\%       & \textbf{1456.7} &-      & \textbf{1481.4} & - \\
\hline
$i=2$&\eqref{eq:ots_mip}           & 1454.8 &0.43\%       & 1456.7 &0.68\%      & 1448.8 & 0.02\%     \\
$i=2$&\eqref{eq:OTS_NP_heur}       & \textbf{1453.6} &-       & 1456.7 &-      & \textbf{1448.6} & - \\
\hline
$i=3$&\eqref{eq:ots_mip}           & 1450.5 & 0.14\%      & 1448.7 &0.01\%      & 1448.6 & -     \\
$i=3$&\eqref{eq:OTS_NP_heur}       & \textbf{1450.4} &-       & \textbf{1448.6} &-      & &           \\
\hline
$i=4$&\eqref{eq:ots_mip}           & 1450.3 & 0.12\%      & 1448.6 & -         & &          \\
$i=4$&\eqref{eq:OTS_NP_heur}       & 1450.3 &-       & &               & &                  \\
\hline
$i=5$&\eqref{eq:ots_mip}           & 1448.8 &0.02\%       & &               & &                   \\
$i=5$&\eqref{eq:OTS_NP_heur}       & 1448.8 &-       & &               & &                      \\
\hline
Final&\eqref{eq:ots_mip}           & 1448.6 &-           & &               & &                \\
      \hline
    \end{tabular}
    \caption{Proposed iterative procedure for a demand instance.}
    \label{tab:iterative_procedure}
\end{table}

Next, we assess the computational efficiency of the three iterative approaches (BN-IT, SR-IT, and LP-IT) by evaluating their performance across 100 instances of the OTS problem. The corresponding computational results are summarized in Table \ref{tab:summary_results}. A comparison between Tables \ref{tab:summary_results} and \ref{tab:summary_results_ss} confirms that the proposed iterative methodology significantly enhances computational performance, regardless of the big-M values used. Specifically, the iterative approach nearly halves the number of unsolved instances for the three approaches and reduces the average computational time by 34.8\%-47.1\% compared to their single-step counterpart approaches. Moreover, the combination of the tightening strategy to refine big-M values using graph simplifications with the iterative methodology, which effectively computes high-quality feasible solutions, yields the best overall computational results.

\begin{figure}
\centering
\begin{tikzpicture}[scale=0.45,font=\LARGE]
	\begin{axis}[	
    width=15cm,
    height=12cm,
    xmin = 0,
    xmax = 3590,
    ymax = 100,
    ymin = 0,
    legend pos = south east,
	legend style={legend cell align=left},	
	clip marker paths=true,	
	xlabel = Time (s),
	ylabel = \# Problems solved]	
    \addplot[line width=2pt,draw=red,dashed] table [x=BN_SS, y=x, col sep=comma] {results_SS_IT_tabla2.csv}; \addlegendentry{\text{BN-SS}}
    \addplot[line width=2pt,draw=blue,dashed] table [x=SR_SS, y=x, col sep=comma] {results_SS_IT_tabla2.csv}; \addlegendentry{\text{SR-SS}}
    \addplot[line width=2pt,draw=green,dashed] table [x=LP_SS, y=x, col sep=comma] {results_SS_IT_tabla2.csv}; \addlegendentry{\text{LP-SS}}
    \addplot[line width=2pt,draw=red] table [x=BN_IT, y=x, col sep=comma] {results_SS_IT_tabla2.csv}; \addlegendentry{\text{BN-IT}}
    \addplot[line width=2pt,draw=blue] table [x=SR_IT, y=x, col sep=comma] {results_SS_IT_tabla2.csv}; \addlegendentry{\text{SR-IT}}
    \addplot[line width=2pt,draw=green] table [x=LP_IT, y=x, col sep=comma] {results_SS_IT_tabla2.csv}; \addlegendentry{\text{LP-IT}}
	\end{axis}	
\end{tikzpicture} 
\caption{Comparison of the number of problems solved over time.}
\label{fig:final_results}
\end{figure}
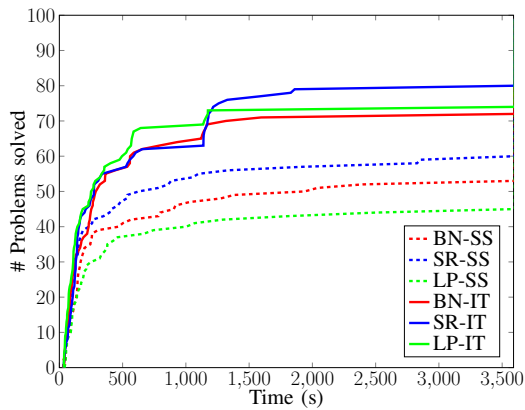

To conclude this numerical analysis, Figure \ref{fig:final_results} illustrates this comparison: the horizontal axis represents computational time, capped at a maximum value of one hour, while the vertical axis shows the cumulative number of instances solved within each time. The results reveal several important observations. First, none of the approaches is able to solve all 100 instances within the one-hour time frame, underscoring the complexity of the problem. Second, for all three sets of big-M values considered, directly solving the OTS model \eqref{eq:ots_mip} using Gurobi is substantially slower than the iterative approach proposed in this study. Finally, the combination of reduced big-M values computed as described in Section \ref{sec:tighter_bigm} with the iterative methodology presented in Section \ref{sec:iterative} shows better computational performance compared to the other approaches, demonstrating its effectiveness in addressing the challenges posed by the OTS problem.

\begin{table}[]
    \centering
    \begin{tabular}{ccccc}
    \hline
    Method  & Av. time (s) & Max Gap & Av. Gap & \# unsolved \\
    \hline
    BN-IT & 1234.16   & 0.43\%     & 0.03\%      & 28  \\
    SR-IT & 1054.95   & 0.39\%     & 0.02\%      & 20  \\
    LP-IT & 1128.42   & 0.43\%     & 0.02\%      & 26  \\
    \hline
    \end{tabular}
    \caption{Summary of computational results for iterative approaches.}
    \label{tab:summary_results}
\end{table}

\section{Conclusions}\label{sec:conclusions}

With the increasing integration of renewable energy sources, network topology modifications have become an essential tool to adapt power systems to highly variable operating conditions. The optimal transmission switching (OTS) problem seeks to identify the set of transmission lines to switch off in order to minimize the operating costs. This problem is formulated as a mixed-integer optimization model that requires sufficiently large big-M constants. However, overly large big-M values can lead to poor relaxations, increased computational burden, and potential numerical issues.

When the network includes a spanning tree of pre-established connected lines to ensure connectivity, appropriate big-M constants can be efficiently determined by solving shortest-path problems. Conversely, when all lines are switchable, computing these constants requires solving the NP-hard longest path problem, making the process more computationally intensive. In this work, we propose a two-fold strategy to address these challenges. First, we leverage graph theory to simplify the network and compute relatively tight big-M values by solving a relaxed version of the longest path problem, which provides a practical upper bound. Second, we introduce an iterative algorithm that incorporates a heuristic procedure to efficiently find high-quality feasible solutions. These solutions are then used to warm-start the integer optimization solver, significantly accelerating convergence.

The proposed strategy is validated on the standard IEEE 118-bus network using 100 demand scenarios. Simulation results demonstrate the effectiveness of the methodology: the number of unsolved instances within one hour is reduced by 57\%, and the average computational time decreases by 44\% compared to existing methods in the technical literature. These results highlight the potential of the proposed approach to enhance the practical applicability of OTS in power systems with high shares of renewable energy.

\bibliographystyle{IEEEtran}
\bibliography{references}

\end{document}